\documentclass[a4paper,12pt]{article}
\begin{document}

\title{Smashed products of topological left quasigroups.}
\author{Sergey V. Ludkowski}
\date{9 October 2020}
\maketitle
\begin{abstract}
In this article smashed products of topological left quasigroups are
scrutinized. Quotients of topological fan quasigroups are
investigated. Skew smashed products of fan quasigroups and their
structure also are studied. \footnote{key words and phrases:
topological
quasigroup; $T_1$; nonlocally compact; product; smashed   }

\end{abstract}
\par Address: Dep. Appl. Mathematics, MIREA - Russian Technological University,
av. Vernadsky 78, Moscow 119454, Russia; \par e-mail:
sludkowski@mail.ru

\section{Introduction.}
\par Topological quasigroups play an important role in algebra, topology,
noncommutative geometry, mathematical physics
\cite{dastqqm80,dastqbcms75,sabininb99,smithb}. They appear
naturally from other areas of mathematics and its applications such
as nonassociative algebra, quantum field theory, topological algebra
\cite{boloktodb,castdoyfior,fernetalmanmath2001,guetze,othmanmarin2017}.
Quasigroups are nonassociative generalizations of groups
\cite{smithb,bruckb,pickert,razm}. A remarkable fact was proved in
the 20-th century that the nontrivial geometry exists if and only if
there exists the corresponding quasigroup.
\par In the article \cite{ludkeimtqgta20} an existence of a left (or
right) invariant nontrivial measure on a topological quasigroup was
investigated. It was proved that if under rather general conditions
the nontrivial left invariant measure exists on the Borel $\sigma
$-algebra of the topological unital quasigroup $G$, then $G$ is
locally compact. \par This article is devoted different types of
products of quasigroups. In Section 2 direct and smashed products of
topological left quasigroups are scrutinized. Quotients of
topological fan quasigroups are investigated. Skew smashed products
of fan quasigroups and their structure also are studied. Examples
are given. It is proved that in particular smashed products of
nontrivial topological groups give topological quasigroups which may
be nonlocally compact. Basic facts on topological quasigroups are
given in Appendix. The obtained results permit to construct ample
classes of topological quasigroups.
\par All main results of this paper are obtained for the first time.

\section{Products of topological quasigroups.}

\par {\bf Definition 1.} Let $G$ be a left (or right) quasigroup. We put:
\par $Com (G) := \{ a\in G: ~ \forall b\in G, ~ ab=ba \} $; $\quad
(1)$
\par $N_l(G) := \{a\in G: ~ \forall b\in G, ~ \forall c\in G, ~ (ab)c=a(bc) \}
$; $\quad (2)$
\par $N_m(G) := \{a\in G: ~ \forall b\in G, ~ \forall c\in G, ~ (ba)c=b(ac)
\} $; $\quad (3)$
\par $N_r(G) := \{a\in G: ~ \forall b\in G, ~ \forall c\in G, ~ (bc)a=b(ca)
\} $; $\quad (4)$
\par $N(G) := N_l(G)\cap N_m(G)\cap N_r(G)$; $\quad (5)$
\par ${\cal C}(G) := Com (G)\cap N(G)$. $\quad (6)$
\par Then $N(G)$ is called a nucleus of $G$
and ${\cal C}(G)$ is called the center of $G$.

\par {\bf Theorem 1.} {\it Let $(G_j, \tau _j)$ be a family of
topological left quasigroups, where $G_j$ is the left quasigroup,
$\tau _j$ is the topology on $G_j$, $j\in J$, $J$ is a set. Then
their direct product $G=\prod_{j\in J}G_j$ relative to the Tychonoff
product topology $\tau _G$ is a topological left quasigroup and
\par  $N(G)=\prod_{j\in J}N(G_j)$ $\quad (7)$ \par and
${\cal C}(G)=\prod_{j\in J}{\cal C}(G_j)$.  $\quad (8)$
\par Moreover,  if $(G_j, \tau _j)$ is $T_1$ for each $j\in J$, then
$(G, \tau )$ is $T_1$, $N(G)$ and ${\cal C}(G)$ are closed in $G$.}
\par {\bf Proof.} The direct product of left quasigroups
satisfies condition $(67)$, hence $G$ is the left quasigroup.
Multiplication and $Div_l$ are jointly continuous relative to the
Tychonoff product topology, consequently, $G$ is the topological
left quasigroup.
\par Note that each element $a\in G$ can be written as $a= \{ a_j: ~ \forall j\in J, ~
a_j\in G_j \}$. From Formulas $(1)$-$(4)$ in Definition 1 we deduce
that
\par  $\quad Com (G) := \{ a\in G: ~ \forall b\in G, ~ ab=ba \}
=\prod_{j\in J} Com (G_j)$, $\quad (9)$ \\ since each $a$ and $b$ in
$G$ have the forms $a= \{ a_j: \forall j\in J, a_j\in G_j \}$ and
$b= \{ b_j: ~ \forall j\in J, b_j\in G_j \} $, and since $ab=ba$ if
and only if
\par $ a_jb_j=b_ja_j$ for all $j\in J$;
\par $\quad N_l(G) := \{a\in G: ~ \forall b\in G, ~ \forall c\in G, ~ (ab)c=a(bc)
\}$\par $\quad = \prod_{j\in J} N_l(G_j)$, $\quad (10)$ \\ since
$(ab)c=a(bc)$ if and only if $(a_jb_j)c_j=a_j(b_jc_j)$ for all $j\in
J$, where $c\in G$ has the form $c= \{ c_j: \forall j\in J, c_j\in
G_j \}; $ and analogously \par  $\quad N_m(G)=\prod_{j\in J}
N_m(G_j)$ $\quad (11)$ \par and  $\quad N_r(G)=\prod_{j\in J}
N_r(G_j)$. $\quad (12)$ \par Thus formulas $(7)$ and $(8)$ follow
from $(9)$-$(12)$, $(5)$ and $(6)$. \par If $(G_j, \tau _j)$ is
$T_1$ for each $j\in J$, then $(G, \tau )$ is $T_1$, since according
to Theorem 2.3.11 in \cite{eng} a product of $T_1$ spaces is a $T_1$
space. From the joint continuity of multiplication it follows that
$Com (G_j)$, $N_l(G_j)$, $N_m(G_j)$, $N_r(G_j)$ are closed in $(G_j,
\tau _j)$, hence $N(G_j)$ and ${\cal C}(G_j)$ are closed in $(G_j,
\tau _j)$ for each $j\in J$, consequently, $N(G)$ and ${\cal C}(G)$
are closed in $(G,\tau )$.

\par {\bf Theorem 2.} {\it Let $(A, \tau _A)$ and $(B, \tau _B)$
be left topological quasigroups, let also $\xi _i: A\times B\times
A\to B$ and $A\times B\ni (a,b)\mapsto \phi _j(a)b\in B$ be
(jointly) continuous univalent maps for each $i\in \{ 1, 2 \} $ and
$j\in \{ 1, 2, 3 \} $, let $\mu : (A\times B)^2\to A\times B$ be a
map such that \par $\mu ((a_1,b_1), (a_2,b_2))= ((a_1,a_2), [(\xi
_1(a_1,b_1,a_2)b_1^{(a_2)}) \xi _2(a_1,b_1,a_2)]^{ \{ a_1 \}
}b_2^{a_1})$ $\quad (13)$ \\ for each $a_1$, $a_2$ in $A$; $b_1$,
$b_2$ in $B$, where \par $b_2^{a_1}:=\phi _1(a_1)b_2$, $\quad (14)$
\par $b_1^{(a_2)}:= \phi _2(a_2)b_1$, $\quad (15)$
\par $b_2^{ \{ a_1 \} }:=\phi _3(a_1)b_2$, $\quad (16)$
\par $\phi _j: A\to {\cal A}(B)$, where ${\cal A}(B)$ denotes the
family of all homeomorphisms from $B$ onto $B$.\par  Then the
cartesian product $C=A\times B$ supplied with the Tychonoff product
topology $\tau _C$ and the map $\mu $ is a topological left
quasigroup $(C,\tau _C)$.}
\par {\bf Proof.} Multiplications in $A$ and $B$ are univalent,
consequently, $\mu $ is univalent and provides multiplication in
$C$. Since $\xi _i(a_1,b_1,a_2)$, $\phi _j(a)b$ are (jointly)
continuous for each $i\in \{ 1, 2 \} $ and $j\in \{ 1, 2, 3 \} $,
$(A, \tau _A)$ and $(B, \tau _B)$ are the left topological
quasigroups, then $\mu $ is (jointly) continuous. Then we consider
the equation:
\par $\mu ((a,b), (x,y))=(c,d)$, $\quad (17)$ \\ where $a\in A$, $c\in A$, $b\in B$, $d\in B$
are arbitrary fixed, $x\in A$, $y\in B$ are elements to be expressed
through $a, b, c, d$. The equation $(17)$ is equivalent to the
system:
\par  $ax=c$ $\quad (18)$ \par and
$[(\xi _1(a,b,x)b^{(x)}) \xi _2(a,b,x)]^{ \{ a \} }y^a=d$. $\quad
(19)$ \\ From $(18)$ it follows that
\par $x=a\setminus c$, $\quad (20)$ \\ since $A$ is the left
quasigroup. Then $(19)$ and $(20)$ imply that
\par $[(\xi _1(a,b,a\setminus c)b^{(a\setminus c)})\xi
_2(a,b,a\setminus c)]^{ \{ a \} }y^a=d$; $\quad (21)$
\par consequently, $y^a=z$ with $z=[(\xi _1(a,b,a\setminus
c)b^{(a\setminus c)})\xi _2(a,b,a\setminus c)]^{ \{ a \} }\setminus
d$ $\quad (22)$ \par and hence $y=[\phi _1(a)]^{-1}z$, $\quad (23)$
\\ since $\phi _1(a)\in {\cal A}(B)$ and $[\phi _1(a)]^{-1}\in {\cal
A}(B)$. Thus the equation $(17)$ has a unique solution $(x,y)$ given
by $(20)$, $(23)$. Denoting it
\par $(x,y)=(a,b)\setminus (c,d)$ $\quad (24)$ \\ we get that $C$
is the left quasigroup. From the (joint) continuity of $\mu $, $\xi
_i$, $\phi _j$ for each $i\in \{ 1, 2 \} $ and $j\in \{ 1, 2, 3 \} $
and Formulas $(20)$, $(22)$-$(24)$ it follows, that the map $Div_l$
is (jointly) continuous, where $Div_l((a,b),(c,d))=(a,b)\setminus
(c,d)$. Thus $(C,\tau _C)$ is the topological left quasigroup.

\par {\bf Definition 2.} The topological left quasigroup provided
by Theorem 2 will also be denoted by $C=A\wp ^{\xi _1, \xi _2, \phi
_1, \phi _2, \phi _3}B$ and called a smashed product of topological
left quasigroups (with smashing factors $\xi _1, \xi _2, \phi _1,
\phi _2, \phi _3$).

\par {\bf Remark 1.} Suppose that the conditions of Theorem 2
are satisfied. It is interesting to study another equation:
\par $\mu ((x,y),(a,b))=(c,d)$, $\quad (25)$ \\ where $a$, $c$ in $A$,
$b$, $d$ in $B$ are arbitrarily fixed, $x$ in $A$ and $y$ in $B$ to
be found as a solution if it exists. It is equivalent to the system:
\par $xa=c$ $\quad (26)$ \par and
$[(\xi _1(x,y,a)y^{(a)}) \xi _2(x,y,a)]^{ \{ x \} }b^x=d$. $\quad
(27)$
\par Suppose in addition that $A$ is a quasigroup. Then there exists
a unique $x=c/a$ in $A$ and $(27)$ takes the form:
\par $[(\xi _1(c/a,y,a)y^{(a)}) \xi _2(c/a,y,a)]^{ \{ c/a \}
}b^{c/a}=d$. $\quad (28)$
\par Even if $B$ is a quasigroup, apparently there may be either no
any or many different solutions $y$ of $(28)$ such that
\par $\xi _1(c/a,y,a)y^{(a)}=f/\xi _2(c/a,y,a)$, $\quad (29)$ \\ where
$f=[\phi _3(c/a)]^{-1}(d/b^{c/a})$.
\\ Indeed, there are $A$ and $B$ such that ${\cal A}(B)- Aut_c(B)\ne
\emptyset $ and ${\sf C}(A\times B\times A,B)-Hom_c(A\times B\times
A,B)\ne \emptyset $, where $Aut_c(B)$ denotes the family of all
continuous automorphisms of $B$, $Hom_c(G,B)$ denotes the family of
all continuous homomorphisms from $G$ into $B$, ${\sf C}(G,B)$
denotes the family of all continuous maps from $G$ into $B$, where
$G$ is a topological left quasigroups (or a quasigroup). Thus
$C=A\wp ^{\xi _1, \xi _2, \phi _1, \phi _2, \phi _3}B$ is the left
topological quasigroup,  but generally it need not be a right
quasigroup.

\par {\bf Example 1.} Let $A$ and $B$ be any topological quasigroups such that
\par $card ({\cal A}(B)-Aut_c(B))\ge card
(Aut_c(B))\ge {\sf c}$ $\quad (30)$ \par and
 $card ( {\sf C}(A\times B\times A,B)-Hom_c(A\times B\times
A,B))$\par $\ge card (Hom_c(A\times B\times A,B))\ge {\sf c}$, $\quad (31)$ \\
where ${\sf c}=card ({\bf R})$. Therefore, there exist $\xi _i$ and
$\phi _j$ (see Theorem 2 and Remark 1) such that $\xi _i\in {\sf
C}(A\times B\times A,B)-Hom_c(A\times B\times A,B)$ and $\phi _j(A)$
is not contained in $Aut_c(B)$ for each $i\in \{ 1, 2 \} $ and $j\in
\{ 1, 2, 3 \} $. Hence $\xi _i$ and $\phi _j$ can be chosen such
that $C=A\wp ^{\xi _1, \xi _2, \phi _1, \phi _2, \phi _3}B$ is the
topological left quasigroup (nonassociative), but not the right
quasigroup.

\par {\bf Example 2.} We take the special orthogonal group
$A=SO(n,{\bf R})$ of the Euclidean space ${\bf R}^n$, the special
linear group $B=SL(m,{\bf R})$ of the Euclidean space ${\bf R}^m$,
where $1<n\le m\in {\bf N}$, $A$ and $B$ are supplied with
topologies induced by the operator norm topology. Therefore for them
Conditions $(30)$ and $(31)$ are satisfied. Hence there exist
topological left quasigroups (nonassociative) $C=A\wp ^{\xi _1, \xi
_2, \phi _1, \phi _2, \phi _3}B$ which are not right quasigroups.
The topological left quasigroup $C$ is locally compact and locally
connected, its small inductive dimension is positive $
\frac{(m^2-1)n(n-1)}{2} = ind (C)$.

\par {\bf Example 3.} Let $l_2$ be the separable Hilbert space over
the complex field ${\bf C}$, where ${\bf C}$ is supplied with the
standard multiplicative norm topology. We consider the unitary group
$A=U(l_2)$ and the general linear group $B=GL(l_2)$ of $l_2$, where
$A$ and $B$ are considered in the topologies inherited form the
operator norm topology. For these groups Conditions $(30)$ and
$(31)$ are satisfied. Therefore there exist their smashed products
$C=A\wp ^{\xi _1, \xi _2, \phi _1, \phi _2, \phi _3}B$ which are
topological left quasigroups (nonassociative), but not right
quasigroups. Moreover, $C$ is locally connected and nonlocally
compact, $ind (C)=\infty $.

\par {\bf Example 4.} Assume that ${\bf F}$ is an infinite
nondiscrete spherically complete field supplied with a
multiplicative norm satisfying the strong triangle inequality
$|a+b|\le \max (|a|, |b|)$ for each $a$ and $b$ in ${\bf F}$ (see
\cite{roo}). \par By $c_0(\alpha ,{\bf F})=X$ is denoted a Banach
space consisting of all vectors $x=(x_j: ~ \forall j\in \alpha , x_j
\in {\bf F})$ such that for each $\epsilon >0$ the cardinality $card
(j\in \alpha : ~ |x_j|>\epsilon )<\aleph _0$ and with a norm
$|x|=\sup_{j\in \alpha } |x_j|$, where $\alpha $ is a (nonvoid) set.
We consider the linear isometry group $A=IL(X)$ and the general
linear group $B=GL(X)$ of $X$ supplied with topologies inherited
from the operator norm topology, where $IL(X)= \{ g\in GL(X): ~
\forall x\in X, ~ |g(x)|=|x| \} $. Evidently, $A$ and $B$ satisfy
Conditions $(30)$ and $(31)$, consequently, there exist their
smashed products $C=A\wp ^{\xi _1, \xi _2, \phi _1, \phi _2, \phi
_3}B$ such that $C$ is a topological left quasigroup and such that
$C$ is not a right quasigroup. \par Therefore $C$ is totally
disconnected, since ${\bf F}$, $X$, $A$ and $B$ are totally
disconnected. If $card (\alpha )<\aleph _0$ and ${\bf F}$ is locally
compact, then $C$ is locally compact. If either $card (\alpha )\ge
\aleph _0$ or ${\bf F}$ is not locally compact, then $C$ is not
locally compact.

\par {\bf Corollary 1.} {\it Let the conditions of Theorem 1
be satisfied and let $G_j$ be compact for each $j\in J_0$ and
locally compact for each $j\in J-J_0$, where $J_0\subset J$. If
$J-J_0$ is a finite set, then $G$ is locally compact. If $J-J_0$ is
infinite and $G_j$ is noncompact for each $j\in J-J_0$, then $G$ is
not locally compact.}
\par {\bf Proof.} By Theorem 1 $G$ is the topological left quasigroup.
In view of Theorem 3.3.13 in \cite{eng} $G$ as the topological space
is locally compact if $J-J_0$ is finite, $G$ is not locally compact
if $J-J_0$ is infinite and $G_j$ is noncompact for each $j\in
J-J_0$.
\par {\bf Corollary 2.} {\it Suppose that the conditions of Theorem 2
are satisfied. Then the smashed product $C=A\wp ^{\xi _1, \xi _2,
\phi _1, \phi _2, \phi _3}B$ is locally compact if and only if $A$
and $B$ are locally compact.}
\par {\bf Proof.} The assertion of this Corollary is evident by
Theorems 3.3.13 in \cite{eng} and 2 above.

\par Next topological quasigroups with some mild restrictions are considered.
They are particular cases of topological left quasigroups. For them
also smashed products are investigated.

\par {\bf Definition 3.} We call $G$ a fan quasigroup
if the unital quasigroup satisfies Conditions $(32)$-$(34)$:
\par $(ab)c=t(a,b,c)a(bc)$ $\quad (32)$ \par  and
$(ab)c=a(bc)p(a,b,c)$ $\quad (33)$ \\ for each $a$, $b$ and $c$ in
$G$, where
\par $t(a,b,c)=t_G(a,b,c)\in N(G)$ and $p(a,b,c)=p_G(a,b,c)\in
N(G)$. $\quad (34)$
\par A minimal closed subgroup $N_0(G)$ in the topological fan quasigroup $G$ containing $t(a,b,c)$ and
$p(a,b,c)$ for each $a$, $b$ and $c$ in $G$ will be called a fan of
$G$.

\par {\bf Corollary 3.} {\it Let the conditions of Theorem 1 be
satisfied and let $G_j$ be a fan quasigroup for each $j\in J$. Then
$G$ is a topological fan quasigroup.}
\par {\bf Proof.} In view of Theorem 1 $G$ is the
topological quasigroup.
\par If $a$, $b$ and $c$ are in $G$, then
\par $(ab)c=\{ (a_jb_j)c_j: ~ \forall j \in J, ~ a_j\in G_j, b_j\in
G_j, c_j\in G_j \} $\par $= \{ t_{G_j}(a_j,b_j,c_j) a_j(b_jc_j): ~
\forall j \in J, ~ a_j\in G_j, b_j\in G_j, c_j\in G_j \} $ \par $=
t_G(a,b,c) a(bc)$ $\quad (35)$ \par and \par $(ab)c =
a(bc)p_G(a,b,c) $ $\quad (36)$ \par with
 $t_G(a,b,c) = \{ t_{G_j}(a_j,b_j,c_j): ~ \forall j \in J, ~
a_j\in G_j, b_j\in G_j, c_j\in G_j \} $ $\quad (37)$ \par and
$p_G(a,b,c) = \{ p_{G_j}(a_j,b_j,c_j): ~ \forall j \in J, ~ a_j\in
G_j, b_j\in G_j, c_j\in G_j \} $. $\quad (38)$
\par Formulas $(35)$-$(38)$ imply that Conditions
$(32)$-$(34)$ are satisfied. On the other hand, Conditions
$(2)$-$(5)$ imply that $N(G)$ is a group. Thus $G$ is the
topological fan quasigroup.

\par {\bf Definition 4.} Let a subquasigroup $H$ of
a quasigroup $G$ satisfy conditions $(39)$ and $(40)$:
\par $xH=Hx$ $\quad (39)$  and
\par $(xy)H=x(yH)$ and $(xH)y=x(Hy)$ and $H(xy)=(Hx)y$ $\quad (40)$ \\ for each $x$ and $y$
in $G$. Then $H$ is called normal in $G$.
\par A family of cosets $\{ bH: ~ b\in G \} $ will be denoted by
$G/_cH$. \par This notation is used below in order to distinguish it
from $b/a=Div_r(a,b)$.

\par {\bf Remark 2.} Assume that $A$ and $B$ are two topological fan
quasigroups; $N_0(A)$ and $N_0(B)$ are fans of $A$ and $B$
respectively.
\par Let also $N$ be a topological group such that
\par $N_0(A)\hookrightarrow N\hookrightarrow N(A)$ and
$ ~ N_0(B)\hookrightarrow N\hookrightarrow N(B)$ $\quad (41)$
\par and let $N$ be normal in $A$ and in $B$ (see also Definitions
1 and 4).
\par Using direct products it is always possible to extend either
$A$ or $B$ to get such a case. In particular, either $A$ or $B$ may
be a group. On $A\times B$ an equivalence relation $\Xi $ is
considered such that
\par $(v\gamma ,b)\Xi (v,\gamma b)$ $\quad (42)$
\\ for every $v$ in $A$, $b$ in $B$ and $\gamma $ in $N$.
\par Let $\phi : A\to {\cal A}(B)$ be a
univalent mapping, $\quad (43)$ \\ where ${\cal A}(B)$ denotes the
family of all homeomorphisms from $B$ onto $B$. If $a\in A$ and
$b\in B$, then we will write $b^a$ for $\phi (a)b$, where $\phi (a)
: B\to B$. Let also \par $\eta _{\phi }: A\times A\times B\to N$,
$\kappa _{\phi } : A\times B\times B\to N$
\par and $\xi
_{\phi } : (A\times B)\times (A\times B)\to N$ \\
be univalent maps written shortly as $\eta $, $\kappa $ and $\xi $
correspondingly such that
\par $(b^u)^v=b^{vu}\eta (v,u,b)$, $~ {\gamma }^u=\gamma $, $~b^{\gamma
}=b$; $\quad (44)$
\par $\eta (v,u,(\gamma _1b)\gamma _2)=\eta (v,u,b)$; $\quad (45)$
\par if $\gamma \in \{ v, u, b \} $ then $\eta (v,u,b)=e$;
\par $(cb)^u=c^ub^u\kappa (u,c,b)$; $\quad (46)$
\par $\kappa (u,(\gamma _1c)\gamma _2,(\gamma _3b)\gamma _4)=\kappa
(u,c,b)$ $\quad (47)$ \par  and  if $\gamma \in \{ u, c, b)$ then
$\kappa (u,c ,b)=e$;
\par $\xi (((\gamma u)\gamma _1,(\gamma _2c)\gamma _3),((\gamma _4v)\gamma _5,
(\gamma _6b)\gamma _7))= \xi ((u,c),(v,b))$ and
\par $\xi ((e,e), (v,b))=e$ and $\xi ((u,c),(e,e))=e$ $\quad (48)$
\\ for every $u$ and $v$ in $A$, $b$, $c$ in $B$,
$\gamma $, $\gamma _1$,...,$\gamma _7$ in $N$, where $e$ denotes the
neutral element in $N$ and in $A$ and $B$.
\par We put
\par $(a_1,b_1)(a_2,b_2)=(a_1a_2,b_2^{a_1}b_1\xi
((a_1,b_1),(a_2,b_2)))$ $\quad (49)$ \\ for each $a_1$, $a_2$ in
$A$, $b_1$ and $b_2$ in $B$. \par The Cartesian product $A\times B$
supplied with such a binary operation $(49)$ will be denoted by
$A\clubsuit ^{\phi , \eta , \kappa , \xi }B$.

\par {\bf Theorem 3.} {\it Let the conditions of Remark 2 be
fulfilled. Then the Cartesian product $A\times B=G$ supplied with
the binary operation $(49)$ is a fan quasigroup.}
\par {\bf Proof.} From the conditions of Remark 2 it follows
that the binary operation $(49)$ is univalent. The group $N$ is
normal in the quasigroups $A$ and $B$ by Conditions $(41)$. Hence
for each $a\in A$ and $\beta \in N$ we get $(a\beta )/a \in N$ and
$a\setminus (\beta a)\in N$, since $aN=Na$ for each $a\in A$.
Analogously for $B$. Thus there are univalent maps
\par $r_{A,a}(\beta )=(a\beta )/a$, $ ~ \check{r}_{A,a}(\beta
)=a\setminus (\beta a)$,
\par $r_{B,b}(\beta )=(b\beta )/b$, $ ~ \check{r}_{B,b}(\beta
)=b\setminus (\beta b)$, \par $r_{A,a}: N\to N$, $ ~
\check{r}_{A,a}: N\to N$, $~ r_{B,b}: N\to N$, $ ~ \check{r}_{B,b}:
N\to N$ \\ for each $a\in A$ and $b\in B$. Evidently \par $r_{A,a}(
\check{r}_{A,a}(\beta ))=\beta $ and $\check{r}_{A,a}(r_{A,a}(\beta
))=\beta $ \\ for each $a\in A$ and $\beta \in N$, and similarly for
$B$.
\par Let $I_1=((a_1,b_1)(a_2,b_2))(a_3,b_3)$ and $I_2=
(a_1,b_1)((a_2,b_2)(a_3,b_3))$, where $a_1$, $a_2$, $a_3$ belong to
$A$; $b_1$, $b_2$, $b_3$ belong to $B$. Then using $(44)$-$(48)$ we
infer that
\par $I_1= ((a_1a_2)a_3,b_3^{a_1a_2}(b_2^{a_1}b_1)\xi ((a_1,b_1),(a_2,b_2))
\xi ((a_1a_2,b_2^{a_1}b_1),(a_3,b_3)))$ and
\par $I_2= (a_1(a_2a_3), (b_3^{a_1a_2}\eta
(a_1,a_2,b_3)b_2^{a_1})$\par $ \kappa (a_1,b_3^{a_2},b_2) \xi
((a_2,b_2),(a_3,b_3))b_1\xi ((a_1,b_1), (a_2a_3, b_3^{a_2}b_2))).$
Therefore \par $I_1=(a,b\alpha )$ with $a=a_1(a_2a_3)$,
$b=(b_3^{a_1a_2}b_2^{a_1})b_1$, \par $\alpha
=\check{r}_{B,b}(p_A(a_1,a_2,a_3))
[p_B(b_3^{a_1a_2},b_2^{a_1},b_1)]^{-1}$\par $ \xi
((a_1,b_1),(a_2,b_2)) \xi ((a_1a_2,b_2^{a_1}b_1),(a_3,b_3))$ and
\par $I_2=(a,b\beta )$ with $\beta =\check{r}_{B,b_1}(\gamma )\xi
((a_1,b_1), (a_2a_3, b_3^{a_2}b_2))$, \\ where $\gamma
=\check{r}_{B,b_2^{a_1}}(\eta (a_1,a_2,b_3))\kappa
(a_1,b_3^{a_2},b_2) \xi ((a_2,b_2),(a_3,b_3))$. Hence
\par  $I_1=I_2p_G$ with $p_G=p_G((a_1,b_1),(a_2,b_2),(a_3,b_3))$ $\quad (50)$
\par and  $I_1=t_GI_2$ with $t_G=t_G((a_1,b_1),(a_2,b_2),(a_3,b_3))$;
\par  $p_G=\beta ^{-1}\alpha $ and $t_G=r_{A,a}(r_{B,b}(p))$, $\quad (51)$
\\ where $G=A\clubsuit ^{\phi , \eta , \kappa , \xi }B$.
\\ Apparently
$t_G ((a_1,b_1),(a_2,b_2),(a_3,b_3))\in N$ and \\ $p_G
((a_1,b_1),(a_2,b_2),(a_3,b_3))\in N$ for each $a_j\in A$, $b_j\in
B$, $j\in \{ 1, 2, 3 \} $, since $\alpha $ and $\beta $ belong to
the group $N$.

\par If $\gamma \in N$ and either $(\gamma ,e)$ or $(e,\gamma )$ belongs
to $\{ (a_1,b_1),(a_2,b_2),(a_3,b_3) \} $, then from the conditions
in Remark 2 and Formulas $(50)$ and $(51)$ it follows that
\par $p_G((a_1,b_1),(a_2,b_2),(a_3,b_3))=e$ and
\par  $t_G((a_1,b_1),(a_2,b_2),(a_3,b_3))=e$, \\ consequently, $(N,e)\cup
(e,N)\subset N(G)$, hence $(N,e)(e,N)=(N,N)\subset N(G)$ by $(42)$.

\par Evidently $(69)$ follows from $(48)$ and $(49)$.
\par Next we consider the following equation
\par $(x,y)(a,b)=(e,e)$, $\quad (52)$ \\ where $a\in A$, $b\in B$ are arbitrarily fixed,
$x\in A$ and $y\in B$ to be calculated.
\par From $(68)$ for the fan quasigroups $A$ and $B$,
$(48)$ and $(49)$ it follows that \par  $x=e/a$, $\quad (53)$
\\ consequently, $b^{(e/a)}y \xi ((e/a,y),(a,b))=e$
and hence \par $y=[b^{(e/a)}\setminus e]/\xi
((e/a,b^{(e/a)}\setminus e),(a,b))]$. $\quad (54)$ \\ Thus $x\in A$
and $y\in B$ given by $(53)$ and $(54)$ provide a unique solution of
$(52)$.
\par Similarly from the following equation
\par $(a,b)(v,z)=(e,e)$, $\quad (55)$ \\ where $a\in A$, $b\in B$
are arbitrarily fixed, $v\in A$ and $z\in B$ are to be found, we
infer that
\par $v=a\setminus e$ $\quad (56)$, \\ consequently,
$z^ab\xi ((a,b),(a\setminus e,z))=e$ and hence \par $z^a= [\xi
((a,b),(a\setminus e,z))]^{-1}/b$ \\ by Conditions $(67)$, $(68)$
and $(43)$ for the fan quasigroups $A$ and $B$. Notice that
$(z^a)^{e/a}=z\eta (e/a,a,z)$, hence by Lemmas 1, 2 and the
conditions of Remark 2
\par $z= \{ [\xi ((a,b),(a\setminus
e,(e/b)^{e/a}))]^{-1}/b \} ^{e/a}/\eta (e/a,a,(e/b)^{e/a}) $. $\quad
(57)$
\\ Thus Formulas $(56)$ and $(57)$ provide a unique solution of
$(55)$.
\par Then it is natural to put $(x,y)=(e,e)/(a,b)$ and
$(v,z)=(a,b)\setminus (e,e)$ and
\par $(a,b)\setminus (c,d)=((a,b)\setminus (e,e))(c,d)
p_G((a,b),(a,b)\setminus (e,e), (c,d))$; $\quad (58)$
\par
$(c,d)/(a,b)=[t_G((c,d),(e,e)/(a,b),(a,b))]^{-1}(c,d)((e,e)/(a,b))$
$\quad (59)$
\par and $e_G=(e,e)$,
where $G=A\clubsuit ^{\phi , \eta , \kappa , \xi }B$. Therefore
Properties $(67)$-$(69)$ and $(32)$-$(34)$ are fulfilled for $G$.

\par {\bf Definition 5.} We call the fan quasigroup
$A\clubsuit ^{\phi , \eta , \kappa , \xi }B$ provided by Theorem 3 a
skew smashed product of the fan quasigroups $A$ and $B$ with
smashing factors $\phi $, $\eta $, $\kappa $ and $\xi $.

\par {\bf Theorem  4.} {\it If $G$ is a $T_1$ topological fan quasigroup,
then its fan $N_0=N_0(G)$ is a normal subgroup. Moreover, if
$N_0\subseteq N_1\subseteq N(G)$, $N_1$ is a closed in $G$ subgroup
satisfying $(39)$, then its quotient $G/_c N_1$ is a $T_1\cap
T_{3.5}$ topological group.}
\par {\bf Proof.} Let $\tau $ be a $T_1$ topology on $G$ relative to
which $G$ is a topological quasigroup. Then each point $x$ in $G$ is
closed, since $G$ is the $T_1$ topological space (see Section 1.5 in
\cite{eng}). From the joint continuity of the multiplication and the
mappings $Div_l$ and $Div_r$ it follows that the nucleus $N(G)$ is
closed in $G$. Therefore the subgroup $N_0$ is the closure of a
subgroup $N_{0,0}(G)$ in $N(G)$ generated by elements $t_G(a,b,c)$
and $p_G(a,b,c)$ for all $a$, $b$ and $c$ in $G$ (see Definitions 1
and 6). According to $(2)$-$(5)$ one gets that $N(G)$ and hence
$N_0$ are subgroups in $G$ satisfying Conditions $(40)$, because
$N_0\subseteq N(G)$.
\par Let $a$ and $b$ belong to $N(G)$ and $x\in G$.
Then $x(x\setminus (ab))=ab$ and \\ $x((x\setminus
a)b)=(x(x\setminus a))b=ab$, consequently, \par $x\setminus
(ab)=(x\setminus a)b$ for each $a$ and $b$ in $N(G)$, $x\in G$. $\quad (60)$ \\
Similarly it is deduced \par $(ab)/x=a(b/x)$ for each $a$ and $b$ in
$N(G)$, $x\in G$. $\quad (61)$ \par Therefore from $(34)$ and $(80)$
and $(60)$ it follows that
\par $((x\setminus a)x)((x\setminus b)x)=(x\setminus
a)(x((x\setminus b)x))p(x\setminus a,x,(x\setminus b)x)$\par
$=(x\setminus (ab))x[p(x,x\setminus b,x)]^{-1}p(x\setminus
a,x,(x\setminus b)x)$,
\\ since $(x\setminus a)(bx)=((x\setminus a)b)x=(x\setminus (ab))x$.
Thus \par $(x\setminus (ab))x= ((x\setminus a)x)((x\setminus
b)x)[p(x\setminus a,x,(x\setminus b)x)]^{-1} p(x,x\setminus b,x)$
for each $a$ and $b$ in $N(G)$, $x\in G$. $\quad (62)$
\par From Identities $(73)$ and $(74)$ it follows that
\par $x\setminus ((u\setminus v)y)=((ux)\setminus (vy))p(u,x,(ux)\setminus (vy)) [p(u,u\setminus v,x)]^{-1}$
$\quad (63)$ \\  for each $u$, $v$, $x$ and $y$ in $G$, since
\par $x\setminus ((u\setminus v)y)=x\setminus (u\setminus (vy))[p(u,u\setminus
v,y)]^{-1}$. \par In particular for $u=a(bc)$ and $v=(ab)c$ with any
$a$, $b$ and $c$ in $G$ we infer using $(34)$ that
$ux=(a(b(cx)))p(b,c,x)p(a,bc,x)$ and $vx=(ab)(cx)p(ab,c,x)$, hence
from $(63)$ and $(93)$ it follows that \par $x\setminus
(p(a,b,c)x)=[p(b,c,x)p(a,bc,x)]^{-1}p(a,b,cx)p(u,x,(ux)\setminus
(vx))$, $\quad (64)$ since \par $x\setminus
(p(a,b,c)x)=[(a(b(cx)))p(b,c,x)p(a,bc,x)]\setminus
[(ab)(cx)p(ab,c,x)]$\par $p(u,x,(ux)\setminus (vx))[p(u,u\setminus
v,x)]^{-1}$, \\ because $u\setminus v=p(a,b,c)\in N(G)$ and
$p(u,u\setminus v,x)=e$.
\par Notice that $(67)$, $(68)$ and $(32)$-$(34)$ imply
$u\setminus (tu)=p$, where $t=t(a,b,c)$, $p=p(a,b,c)$, $u=a(bc)$ for
any $a$, $b$ and $c$ in $G$. Let $z\in G$, then there exists $x\in
G$ such that $z=ux$, that is $x=u\setminus z$. Therefore we deduce
that
\par  $z\setminus (tz)=[x\setminus (px)]p(u,u\setminus
(tu),x)[p(u,x,(ux)\setminus (tux))]^{-1}$, $\quad (65)$ \\
since $t\in N(G)$, $p\in N(G)$, $(u\setminus (tu))x=(u\setminus
(tux))[p(u,u\setminus (tu),x)]^{-1}$ by $(74)$; $~x\setminus
(u\setminus (tux))=[(ux)\setminus (tux))]p(u,x,(ux)\setminus (tux))$
by $(73)$. Thus from Identities $(62)$, $(64)$ and $(65)$ it follows
that a group $N_{0,0}=N_{0,0}(G)$ generated by $ \{ p(a,b,c), ~
t(a,b,c): a\in G, ~ b\in G, ~ c\in G \} $ satisfies Condition
$(39)$. From the joint continuity of multiplication and the mappings
$Div_l$ and $Div_r$ it follows that the closure $N_0$ of $N_{0,0}$
also satisfies $(39)$. Thus $N_0$ is a closed normal subgroup in
$G$.
\par Since $N_1$ satisfies $(39)$ and $N_1$ is the subgroup in $G$
such that $N_0\subseteq N_1\subseteq N(G)$, then $N_1$ is normal in
$G$. Therefore a quotient quasigroup $G/_cN_1$ exists consisting of
all cosets $aN_1$, where $a\in G$.
\par  Then from
Conditions $(34)$, $(39)$ and $(40)$ it follows that for each $a$,
$b$, $c$ in $G$ the identities take place \par
$(aN_1)(bN_1)=(ab)N_1$ and
\par $((aN_1)(bN_1))(cN_1)=(aN_1)((bN_1)(cN_1))$ and
$eN_1=N_1$, \\ since $p_G(a,b,c)\in N_0\subseteq N_1$ and
$t_G(a,b,c)\in N_0\subseteq N_1$ for all $a$, $b$ and $c$ in $G$.
\par In view of Lemmas 1 and 2 $(aN_1)\setminus e= e/(aN_1)$,
consequently, for each $aN_1\in G/_c N_1$ a unique inverse
$(aN_1)^{-1}$ exists. Thus the quotient $G/_cN_1$ of $G$ by $N_1$ is
a group. Since the topology $\tau $ on $G$ is $T_1$ and $N_1$ is
closed in $G$ by the conditions of this theorem, then the quotient
topology $\tau _q$ on $G/_c N_1$ is also $T_1$. By virtue of Theorem
8.4 in \cite{hew} this implies that $\tau _q$ is a $T_1\cap T_{3.5}$
topology on $G/_c N_1$.

\par {\bf Corollary 4.} {\it Suppose that the conditions of Remark 2
are fulfilled and $A$ and $B$ are topological $T_1$ fan quasigroups
and smashing factors $\phi $, $\eta $, $\kappa $, $\xi $ are jointly
continuous by their variables. Suppose also that $A\clubsuit ^{\phi
, \eta , \kappa , \xi }B$ is supplied with a topology induced from
the Tychonoff product topology on $A\times B$. Then $A\clubsuit
^{\phi , \eta , \kappa , \xi }B$ is a topological $T_1$ fan
quasigroup.}

\par {\bf Remark 3.} In particular, it is possible to consider a topological quasigroup
$G$ satisfying the condition:
\par there exists a compact
subgroup $N_0=N_0(G)$ in $N(G)$ such that \par $t_G(a,b,c)\in N_0$
and $p_G(a,b,c)\in N_0$ for every $a$, $b$ and $c$ in $G$. $\quad
(66)$

\par {\bf Corollary 5.} {\it If the conditions of Corollary 4
are satisfied and quasigroups $A$ and $B$ are locally compact, then
$A\clubsuit ^{\phi , \eta , \kappa , \xi }B$ is locally compact.
Moreover, if $A$ and $B$ satisfy Condition $(66)$ and ranges of
$\eta $, $\kappa $, $\xi $ are contained in $N_0(A)N_0(B)$, then
$A\clubsuit ^{\phi , \eta , \kappa , \xi }B$ satisfies Condition
$(66)$.}
\par {\bf Proof.} Corollaries 4 and 5 follow immediately from
Theorems 2.3.11, 3.2.4, 3.3.13 in \cite{eng}, Lemma 2.6 in
\cite{ludkeimtqgta20} and Theorems 3, 4, Lemma 3, since
\par $N_0(A)N_0(B)\subseteq N\subseteq N(A)\cap N(B)$ and because
$N_0(A)N_0(B)$ is a compact subgroup in $A\clubsuit ^{\phi , \eta ,
\kappa , \xi }B$.

\par {\bf Remark 4.} From Theorems 1, 3, 4 and Corollaries
4, 5 it follows that taking nontrivial $\phi $, $\eta $, $\kappa $
and $\xi $ and starting even from groups with nontrivial $N(G_j)$ or
$N(A)$ and $G_j/_cN(G_j)$ or $A/_cN(A)$ it is possible to construct
new fan quasigroups with nontrivial $N_0(G)$ and ranges $t_G(G,G,G)$
and $p_G(G,G,G)$ of $t_G$ and $p_G$ may be infinite and nondiscrete.
With suitable smashing factors $\phi $, $\eta $, $\kappa $ and $\xi
$ and with nontrivial fan quasigroups or groups $A$ and $B$ it is
easy to get examples of fan quasigroups in which $e/a\ne a\setminus
e$ for an infinite family of elements $a$ in $A\clubsuit ^{\phi ,
\eta , \kappa , \xi }B$. It can be seen that the smashed product of
topological quasigroups is a generalization of a semidirect product
of topological groups.

\par {\bf Conclusion.} The results of this article can be used for
further investigations of structure of topological quasigroups,
their smashed and skew smashed products, homogeneous spaces
associated with quasigroups \cite{fedumn86}, measures on homogeneous
spaces and noncommutative manifolds \cite{pickert}, microbundles
\cite{milnor64,ludkmbtrta19}. It is worth to mention possible
applications in mathematical coding theory, techniques such as
information flows analysis and systems with distributed memory
\cite{blautrctb,magglebrtj519,srwseabm14}. Indeed, codes are
frequently based on topological-algebraic binary systems. Naturally,
they also can be utilized in harmonic analysis on nonassociative
algebras \cite{razm}, quasigroups, their representation theory
\cite{smithb}, geometry, mathematical physics, quantum field theory,
gauge theory, PDEs, etc.
\cite{boloktodb,guetze,gilmurr,girard,guespr,ludkcvee16,ludkrimut13}.

\section{Appendix. Basics on topological quasi-groups.}
\par For convenience we remind a definition,
though a reader familiar with \cite{smithb,bruckb,eng} may skip
Definition 6.
\par {\bf Definition 6.}  Let $G$ be a set with multiplication
(that is a univalent binary operation) $G^2\ni (a,b)\mapsto ab \in
G$ defined on $G$ such that
\par  for each $a$ and $b$ in $G$ there is a unique $x\in
G$ with $ax=b$.  $\quad (67)$
\par A set $G$ possessing multiplication and satisfying condition
$(67)$ is called a left quasigroup. Symmetrically it is considered
\par a unique $y\in G$ exists satisfying $ya=b$. $\quad (68)$
\par Then the set $G$ possessing multiplication and satisfying condition
$(68)$ is called a right quasigroup.
\par The maps in $(67)$ and $(68)$ are denoted by $x=a\setminus b=Div_l(a,b)$ and
$y=b/a=Div_r(a,b)$ correspondingly.
\par If $G$ is a left and right quasigroup, then it is called a quasigroup.
\par Let $\tau $ be a topology on the left (or right) quasigroup $G$ such that
multiplication $G\times G\ni (a,b)\mapsto ab\in G$ and the mapping
$Div_l(a,b)$ (or $Div_r(a,b)$ respectively) are jointly continuous
relative to $\tau $, then $(G, \tau )$ will be called a topological
left (or right respectively) quasigroup. If $G$ is a topological
left and right quasigroup, then it is called a topological
quasigroup.
\par A set $G$ possessing multiplication is called a groupoid.
\par If there exists a neutral (i.e. unit) element $e_G=e\in G$: $~eg=ge=g$
for each $g\in G$, $\quad (69).$ \\ then the groupoid $G$ is called
unital.
\par If $A$ and $B$ are subsets in $G$, then $A-B$ means the difference of them $A-B=\{
a\in A: ~a \notin B \} $.

\par {\bf Lemma 1.} {\it If $G$ is a fan quasigroup, then
for each $a$, $b$ and $c$ in $G$ the following identities are
fulfilled:
\par $ ~ b\setminus e=t(e/b,b,b\setminus e)(e/b)$;  $\quad (70)$
\par $b\setminus e=(e/b)p(e/b,b,b\setminus e)$;  $\quad (71)$
\par $(a\setminus e)b=t(e/a,a,a\setminus e)[t(e/a,a,a\setminus b)]^{-1}(a\setminus
b)$;  $\quad (72)$
\par $(a\setminus b)= (a\setminus e)bp(a,a\setminus e,b)$;
\par $(bc)\setminus a=(c\setminus (b\setminus
a))[p(b,c,(bc)\setminus a)]^{-1}$; $\quad (73)$
\par  $(a\setminus b)c=(a\setminus (bc))[p(a,a\setminus b,c)]^{-1}$;
 $\quad (74)$
\par  $(ab)\setminus e = (b\setminus e)(a\setminus
e)[t(a,b,b\setminus e)]^{-1}t(ab,b\setminus e,a\setminus e)$; $\quad
(75)$
\par $b(e/a)=(b/a)p(b/a,a,a\setminus e)[p(e/a,a,a\setminus e)]^{-1} $;
\par $(b/a)= [t(b,e/a,a)]^{-1}b(e/a)$;  $\quad (76)$
\par  $a/(bc)=t(a/(bc),b,c)((a/c)/b)$;  $\quad (77)$
\par $c(b/a)=t(c,b/a,a)(cb)/a$;  $\quad (78)$
\par $e/(ab)=[p(e/b,e/a,ab)]^{-1}p(e/a,a,b)(e/b)(e/a)$. $\quad (79)$ }
\par {\bf Proof.} Note that $N(G)$ is a subgroup in $G$ due to Conditions
$(2)$-$(5)$. Then Conditions $(67)$-$(69)$ imply that
\par $b(b\setminus a)=a$, $~b\setminus (ba)=a$;  $\quad (80)$
\par $(a/b)b=a$, $~(ab)/b=a$  $\quad (81)$ \\
for each $a$ and $b$ in any quasigroup $G$. Using Conditions
$(32)$-$(34)$ and Identities $(80)$ and $(81)$ we deduce that
\par $e/b=(e/b)(b(b\setminus e))  = [t(e/b,b,b\setminus e)]^{-1}(b\setminus e)
$\\ which leads to $(70)$.
\par Let $c=a\setminus b$, then from Identities $(70)$ and $(80)$ it
follows that \par $(a\setminus e)b=t(e/a,a,a\setminus
e)(e/a)(ac)$\par $=t(e/a,a,a\setminus e)[t(e/a,a,a\setminus b)]^{-1}
((e/a)a)(a\setminus b)$\\
which taking into account $(81)$ provides $(72)$.
\par On the other hand, $b\setminus e=((e/b)b)(b\setminus e) =
(e/b)(b(b\setminus e))p(e/b,b,b\setminus e)$ that gives $(71)$.
\par Let now $d=b/a$, then Identities $(71)$ and $(81)$ imply
that  \par $b(e/a)=(da)(a\setminus e)[p(e/a,a,a\setminus
e)]^{-1}$\par $=(b/a)p(b/a,a,a\setminus e)[p(e/a,a,a\setminus e)]^{-1} $ \\
which demonstrates $(76)$.
\par Next we infer from $(32)$-$(34)$ and $(80)$ that
\par $b(c((bc)\setminus a))=(bc)((bc)\setminus a)[p(b,c,(bc)\setminus
a)]^{-1}= a[p(b,c,(bc)\setminus a)]^{-1}$, hence $c((bc)\setminus
a)=(b\setminus a)[p(b,c,(bc)\setminus a)]^{-1}$ that implies $(73)$.
\par Symmetrically it is deduced that
$(a/(bc))b)c=t(a/(bc),b,c)a$, consequently,
$(a/(bc))b=t(a/(bc),b,c)(a/c)$. From the latter identity it follows
$(77)$. \par Evidently, formulas \par $a((a\setminus
b)c)=(a(a\setminus b))c[p(a,a\setminus b,c)]^{-1}=bc[p(a,a\setminus
b,c)]^{-1}$ and \par $(c(b/a))a=t(c,b/a,a)cb$ \\ imply $(74)$ and
$(78)$ correspondingly.
\par From $(34)$ we infer that \par $(ab)((b\setminus
e)(a\setminus e))=[t(ab,b\setminus e,a\setminus
e)]^{-1}t(a,b,b\setminus e)$, since by $(80)$
\par $(a(b(b\setminus e)))(a\setminus e)=e$. \\ This together with
$(67)$ and $(68)$ implies $(75)$.
\par Analogously form $(32)$ we deduce that
\par $((e/b)(e/a))(ab)=[p(e/a,a,b)]^{-1}p(e/b,e/a,ab)$, since by
$(81)$
\par $(e/b)(((e/a)a)b)=e$.
\\ Finally applying $(67)$ and $(68)$ we get Identity $(79)$.

\par {\bf Lemma 2.} {\it Assume that $G$ is a fan quasigroup.
Then for every $a$, $a_1$, $a_2$, $a_3$ in $G$ and $z_1$, $z_2$,
$z_3$ in ${\cal C}(G)$, $b\in N(G)$:
\par $t(z_1a_1,z_2a_2,z_3a_3)=t(a_1,a_2,a_3)$; $\quad (82)$
\par $p(z_1a_1,z_2a_2,z_3a_3)=p(a_1,a_2,a_3)$; $\quad (83)$
\par  $t(a,a\setminus e,a)a=ap(a,a\setminus e,a)$;  $\quad (84)$
\par $t(a,e/a,a)a=ap(a,e/a,a)$; $\quad (85)$
\par  $p(a,a\setminus e,a)t(e/a,a,a\setminus e)=e$; $\quad (86)$
\par $t(a_1,a_2,a_3b)=t(a_1,a_2,a_3)$; $\quad (87)$
\par $p(ba_1,a_2,a_3)=p(a_1,a_2,a_3)$;  $\quad (88)$
\par $t(ba_1,a_2,a_3)=bt(a_1,a_2,a_3)b^{-1}$; $\quad (89)$
\par $p(a_1,a_2,a_3b)=b^{-1}p(a_1,a_2,a_3)b$. $\quad (90)$ }
\par {\bf Proof.} Since $(a_1a_2)a_3=t(a_1,a_2,a_3)a_1(a_2a_3)$ and
$t(a_1,a_2,a_3)\in N(G)$ for every $a_1$, $a_2$, $a_3$ in $G$, then
\par $t(a_1,a_2,a_3)=((a_1a_2)a_3)/(a_1(a_2a_3))$. $\quad (91)$ \par Therefore, for
every $a_1$, $a_2$, $a_3$ in $G$ and $z_1$, $z_2$, $z_3$ in ${\cal
C}(G)$ we infer that \par $t(z_1a_1,z_2a_2,z_3a_3)=
(((z_1a_1)(z_2a_2))(z_3a_3))/((z_1a_1)((z_2a_2)(z_3a_3)))=$\par $
((z_1z_2z_3)((a_1a_2)a_3))/((z_1z_2z_3)(a_1(a_2a_3)))=
((a_1a_2)a_3)/(a_1(a_2a_3))$, $\quad (92)$ \\ since
$b/(qa)=q^{-1}b/a$ and $b/q=q\setminus b=bq^{-1}$  \\ for each $q\in
{\cal C}(G)$, $ ~ a$ and $b$ in $G$, because ${\cal C}(G)$ is the
commutative group satisfying Conditions $(1)$ and $(6)$. Thus
$t(z_1a_1,z_2a_2,z_3a_3)=t(a_1,a_2,a_3)$.
\par Symmetrically we get \par  $p(a_1,a_2,a_3)=(a_1(a_2a_3))\setminus ((a_1a_2)a_3)$ and
\par $p(z_1a_1,z_2a_2,z_3a_3)=
((z_1a_1)((z_2a_2)(z_3a_3)))\setminus
(((z_1a_1)(z_2a_2))(z_3a_3))$\par
$=((z_1z_2z_3)(a_1(a_2a_3)))\setminus ((z_1z_2z_3)((a_1a_2)a_3))=
(a_1(a_2a_3))\setminus ((a_1a_2)a_3)$ $\quad (93)$ \\ that provides
$(83)$.
\par From Formulas $(91)$ and $(70)$ it follows that \par $t(a,a\setminus
e,a)=((a(a\setminus e))a)/(a((a\setminus
e)a))=a/[at(e/a,a,a\setminus e)]$, consequently, \par
\par  $t(a,a\setminus e,a)at(e/a,a,a\setminus e)=a$. $\quad (94)$ \\ Then
from Formulas $(93)$, $(80)$ and Conditions $(32)$-$(34)$ we deduce
that \par $p(a,a\setminus e,a)=(a((a\setminus e)a))\setminus
((a(a\setminus e))a)= \{ [t(a,a\setminus e,a)]^{-1}a \} \setminus
a$,
\\ which implies $(84)$. Identities $(84)$ and $(94)$ lead to
$(86)$. Next using $(93)$ and $(34)$ we deduce that \par
$p(a,e/a,a)=[a((e/a)a)]\setminus [(a(e/a))a]=a\setminus [t(a,e/a,a)a
]$ \\ that implies $(85)$. From $(34)$ we get that
\par $((a_1a_2)a_3)b=(a_1a_2)(a_3b)=(t(a_1,a_2,a_3b)a_1(a_2a_3))b$, \\
from which and $(81)$ and $(91)$ Identity $(87)$ follows, because
$b\in N(G)$. Then
\par $b((a_1a_2)a_3)=((ba_1)a_2)a_3=b(a_1(a_2a_3)p(ba_1,a_2,a_3))$ \\
and $(80)$ and $(93)$ imply Identity $(88)$. Symmetrically we deduce
\par $b((a_1a_2)a_3)=t(ba_1,a_2,a_3))b(a_1(a_2a_3))$ and
\par $((a_1a_2)a_3)b=(a_1(a_2a_3))bp(a_1,a_2,a_3b)$ \\ that together with
$(91)$ and $(93)$ imply Identities $(89)$ and $(90)$.

\par {\bf Lemma 3.} {\it If $(G, \tau )$ is a topological quasigroup, then
the functions $t(a_1,a_2,a_3)$ and $p(a_1,a_2,a_3)$ are jointly
continuous in $a_1$, $a_2$, $a_3$ in $G$.}
\par {\bf Proof.} This follows immediately from Formulas $(91)$, $(93)$ and
Definitions 1, 6.

\end{document}